\documentclass[12pt]{amsart}

\usepackage{amsfonts,amsmath,amssymb,amsthm}
\usepackage[dvips]{graphicx,color}

\newtheorem{theorem}{Theorem}
\newtheorem{proposition}[theorem]{Proposition}

\newtheorem{corollary}[theorem]{Corollary}

\theoremstyle{definition}

\theoremstyle{remark}

\newcommand{\comp}{\vDash}            
\DeclareMathOperator{\Co}{Co}   

\DeclareMathOperator{\Gtype}{Type}

\newcommand{\integers}{\mathbb{Z}}   
\def\field{\integers}

\newcommand{\GDesAlg}{\mathcal{D}}

\title[Semigroups, wreath products, and descent algebras]{A semigroup approach to wreath-product extensions of Solomon's descent algebras}
\author{Samuel K. Hsiao}


 \date{October 15, 2007}

\address{Mathematics Program, Bard College, Annandale-on-Hudson, NY, 12504}

\begin{document}
\maketitle

\begin{abstract}
There is a well-known combinatorial definition, based on ordered set partitions, of the semigroup of faces of the braid arrangement. We generalize this definition to obtain a semigroup $\Sigma_n^G$ associated with $G\wr S_n$, the wreath product of the symmetric group $S_n$ with an arbitrary group $G$. Techniques of Bidigare and Brown are adapted to construct an anti-homomorphism from the $S_n$-invariant subalgebra of the semigroup algebra of $\Sigma_n^G$ into the group algebra of $G\wr S_n$. The generalized descent algebras of Mantaci and Reutenauer are obtained as homomorphic images when $G$ is abelian.
\end{abstract}

\section{Introduction}
A celebrated result of Solomon \cite{Sol} reveals the existence of an intriguing subalgebra, called the {\it descent algebra}, inside the group algebra of any Coxeter group. Solomon's descent algebras have been studied extensively, particularly for the symmetric group $S_n$, in which case there are rich connections to the study of free Lie algebras \cite{GR, R}, quasisymmetric functions \cite{G, MalR}, shuffling \cite{BD}, and Markov chains on faces of hyperplane arrangements \cite{B}. 

Bidigare \cite{Bid} gave a geometrically motivated proof Solomon's result for $S_n$, which involves analyzing an action of $S_n$ on the semigroup $\Sigma_n$ of faces of the braid arrangement (Coxeter complex of $S_n$). In this paper we generalize Bidigare's approach to $G\wr S_n$, the wreath product of $S_n$ with an arbitrary finite group $G$. We construct a class of semigroups that can be viewed as wreath-product analogs of $\Sigma_n$. Each semigroup $\Sigma_n^G$ is defined in terms of ordered set partitions of $\{1,2,\ldots, n\}$ decorated with elements from $G$, generalizing the combinatorial definition of $\Sigma_n$. Unlike the face semigroup of the braid arrangement (or of a hyperplane arrangement in general), elements of $\Sigma_n^G$ are not necessarily idempotent. Instead, they satisfy the identities
\begin{equation}\label{E:semigroup}
 x^{|G|+1} = x \quad \text{ and } \quad xyx^{|G|} = xy
 \end{equation}
for all $x,y \in \Sigma_n^G$. When $|G|=1$ these identities define left regular bands. In general they define a class of semigroups in which every element belongs to a subgroup of exponent $|G|$ (hence they are ``completely regular") and whose idempotents form a left regular band. 

We introduce an $S_n$-action on the semigroup algebra $\field\Sigma_n^G$, for which the invariant subalgebra $(\field\Sigma_n^G)^{S_n}$ has a basis $(\sigma_\alpha)$ indexed by $G$-compositions (these generalize the notion of ``descent set"). The group algebra $\field[G\wr S_n]$ also contains a $\field$-submodule, defined analogously to Solomon's descent algebra, consisting of a natural basis $(X_\alpha)$ indexed by $G$-compositions.  Our main result is as follows:

\begin{theorem}
The $\field$-module map $f:(\field\Sigma_n^G)^{S_n}\to\field[G\wr S_n]$ given by $f(\sigma_\alpha) = X_\alpha$ is an injective anti-homomorphism of algebras. 
\end{theorem}

It follows that the image of $f$ is a subalgebra of $\field[G\wr S_n]$. For abelian groups $G$ these subalgebras turn out to be the generalized descent algebras of Mantaci and Reutenauer  \cite{MR}. Recent work by Hohlweg \cite{H} together with  Baumann \cite{BH} and Bergeron \cite{BerH} studies these algebras from representation-theoretic and Hopf-algebraic viewpoints. For arbitrary groups $G$ these algebras have been considered by Novelli and Thibon \cite{NT} in their work generalizing the study of free quasisymmetric functions to the context of colored permutations.

Our approach, in addition to providing an elementary and concise proof of the existence of generalized descent algebras, suggests a framework in which to develop wreath-product versions of algebraic, combinatorial, and probabilistic results on left regular bands. For instance, could Brown's theory of Markov chains on left regular bands \cite{B}, Saliola's work on the structure of semigroup algebras of left regular bands \cite{S, S2}, or the Hopf algebraic and combinatorial constructions of Aguiar and Mahajan \cite{AM}, be extended to semigroups axiomatized by \eqref{E:semigroup}?

We learned of Bigidare's (unpublished) result through a paper of Brown \cite[Theorem~7]{B}, in which a geometric version of Bidigare's proof is given. While our proof is purely algebraic, it has Brown's argument at its core.

Throughout this paper we assume that $G$ is a finite group and denote its identity element by $e$. However, it should be noted that the definitions of our main objects of study, namely $\Sigma_n^G, (\Sigma_n^G)^{S_n}, G\wr S_n,$ and $\GDesAlg(G\wr S_n)$, still make sense when $G$ is an infinite group, or even just a semigroup with identity. For such $G$, Theorem~1 still holds, but  Formulas~\eqref{E:semigroup} might fail.

\section{Preliminaries}

\subsection{The semigroup of ordered $G$-partitions}
Fix a positive integer $n$, and let $[n] = \{1,2,\ldots, n\}$. An {\em ordered partition} (also called block partition) of $[n]$ is a tuple $(B_1, \ldots, B_k)$ of nonempty pairwise disjoint sets whose union is $[n]$.

An {\em ordered $G$-partition of $[n]$} is a tuple $((B_1, g_1), \ldots, (B_k, g_k))$ such that  $(B_1, \ldots, B_k)$ is an ordered partition of $[n]$ and $g_i\in G$ for all $i\in [k]$.

Let $\Sigma_n^G$ denote the set of ordered $G$-partitions of $[n]$. Define multiplication in $\Sigma_n^G$ by
\[
 ((B_1, g_1), \ldots, (B_k, g_k)) ((C_1, h_1), \ldots, (C_\ell, h_\ell)) = 
\begin{array}{c} 
( (B_1\cap C_1, h_1 g_1), \ldots, (B_1 \cap C_\ell, h_\ell g_1), \\
 (B_2\cap C_1, h_1 g_2), \ldots, (B_2 \cap C_\ell, h_\ell g_2), \\
\vdots \\
(B_k \cap C_1, h_1 g_k),\ldots, (B_k \cap C_\ell, h_\ell g_k) )
\end{array}
\]
where empty intersections are omitted. This gives $\Sigma_n^G$ the structure of a semigroup (with identity element $(([n], e))$) satisfying Formulas \eqref{E:semigroup}. If $|G| = 1$ then $\Sigma_n^G$ is isomorphic to the face semigroup of the braid arrangement. See \cite{B} for details.

\subsection{The invariant subalgebra}
The action of the symmetric group $S_n$ on $[n]$ induces an action on $\Sigma_n^G$. For example,  $\pi \cdot ((\{1,3\}, g_1), (\{2\}, g_2)) = ((\{\pi(1),\pi(3)\},g_1),(\{\pi(2)\}, g_2))$ for any $\pi \in S_3$. This action extends linearly to the semigroup algebra $\field\Sigma_n^G$. Consider the subalgebra of invariants under the action of $S_n$:
\[(\field\Sigma_n^G)^{S_n} = \{ P \in \field\Sigma_n^G~|~\pi\cdot P = P \text{ for all } \pi \in S_n\}.\]
That $(\field\Sigma_n^G)^{S_n}$ is a subalgebra of $\field\Sigma_n^G$ is a consequence of the observation that $\pi\cdot(P Q) = (\pi\cdot P) (\pi\cdot Q)$ for all $\pi \in S_n$ and $P, Q \in \Sigma_n^G$. 

As a $\integers$-module $(\field\Sigma_n^G)^{S_n}$ is free with a basis indexed by $G$-compositions. By a {\em $G$-composition of $n$} we mean a sequence $\alpha = ((a_1, g_1), \ldots, (a_k, g_k))$ such that $(a_1, \ldots, a_k)$ is a composition of $n$, i.e.\ a list of positive integers summing to $n$, and $g_i \in G$ for all $i\in [k]$. In this case we write $\alpha \comp_G n$ and $\ell(\alpha) = k$. The {\em type} of an ordered $G$-partition is the $G$-composition defined by 
\[\Gtype(((B_1,g_1),\ldots,(B_k,g_k))) = ((|B_1|, g_1), \ldots, (|B_k|,g_k)).\] 
For $\alpha \comp_G n$, let
\[ \sigma_\alpha = \sum_{P \in \Sigma_n^G : \Gtype(P) =\alpha} P.\]
Clearly $(\sigma_\alpha)_{\alpha\comp_G n}$ is a basis for $(\field\Sigma_n^G)^{S_n}$. 

\subsection{Multiplication rule for the invariant subalgebra} For $\alpha,\beta,\gamma \comp_G n$, the coefficient of $\sigma_\gamma$ in the product $\sigma_\alpha \sigma_\beta$ is just the 
number of ways of writing an arbitrary $R\in \Sigma_n^G$ of type $\gamma$ as a product $R = P Q$ where $\Gtype(P) = \alpha$ and $\Gtype(Q) = \beta$. Thus, by the multiplication rule for ordered $G$-partitions, we obtain the following multiplication rule inside $(\field\Sigma_n^G)^{S_n}$. Consider all $k\times l$ matrices $M$ whose entries are of the form $M_{ij} = 0$ or $M_{ij} = (a,g)$ where $a$ is a positive integer and $g \in G$. Let $|0| = 0$ and $|(a,g)| = a$, and call $g$ the {\it color} of $(a,g)$. Say that $M$ is {\it compatible} with $\alpha$ and $\beta$, where $\alpha = ((a_1, g_1),\ldots, (a_k,g_k)) \comp_G n$ and  $\beta = ((b_1, h_1),\ldots, (b_\ell, h_\ell))\comp_G n$, if the following conditions are satisfied:
\begin{itemize}
\item[(a)] For all $i\in [k]$, $\sum_{j=1}^\ell |M_{ij}| = a_i$,
\item[(b)] For all $j\in [\ell]$, $\sum_{i=1}^k |M_{ij}| = b_j$,
\item[(c)] For all $i\in [k]$ and $j\in [\ell]$, if $M_{ij} \ne 0$ then $M_{ij}$ has color $h_j g_i$.
\end{itemize}
For a compatible matrix $M$, let $M'$ denote the $G$-composition obtained by reading the entries of $M$ row-by-row, omitting entries that are $0$. For example, the following matrix is compatible with  $\alpha = ((4,g_1), (6,g_2))$ and $\beta = ((3,h_1), (5,h_2), (2,h_3))$:
\[
M = \left(\begin{matrix}
(2, h_1 g_1) & 0 & (2, h_3 g_1) \\
(1,h_1 g_2) & (5, h_2 g_2) & 0  
\end{matrix}\right)
\]
Here, $M' = ((2,h_1g_1),(2,h_3g_1),(1,h_1g_2),(5,h_2g_2))$.

\begin{proposition}\label{P:mult}
Given $G$-compositions $\alpha = ((a_1, g_1),\ldots, (a_k,g_k))$ and $\beta = ((b_1, h_1),\ldots, (b_\ell, h_\ell))$ of $n$, we have
\[ \sigma_\alpha \sigma_\beta = \sum_M \sigma_{M'}\]
where the sum is over all matrices compatible with $\alpha$ and $\beta$.
\end{proposition}
 
When $G$ is abelian, Proposition~\ref{P:mult} is equivalent to the formula for multiplication inside the generalized descent algebra obtained  by Mantaci and Reutenauer \cite[Corollary~6.8]{MR}. This formula is originally due to  Garsia and Remmel \cite{GRem} for the descent algebra of $S_n$.

\subsection{The $G$-descent algebra}
Consider the right permutation action of $S_n$ on $G^{[n]}$, the group of functions from $[n]$ to $G$ with multiplication given by $(gh)(i) = g(i) h(i)$ for $g,h\in G^{[n]}$ and $i\in [n]$. A permutation $\pi\in S_n$ takes $g \in G^{[n]}$ to $g\cdot \pi$, where $(g \cdot \pi)(i) = g(\pi(i))$. Using this action we construct the wreath product $G \wr S_n$. As a set, $G\wr S_n = S_n \times G^{[n]}$. Its group operation is given by $(\pi, g)*(\tau, h) = (\pi\tau, (g\cdot \tau)h)$. It will be convenient to represent an element $(\pi, g) \in G \wr S_n$ by $((\pi_1, g_1), \ldots, (\pi_n, g_n))$, where $\pi_i = \pi(i)$ and $g_i = g(i)$ for $i\in [n]$. With this notation,
\begin{equation}
\label{E:desmult} ((\pi_1, g_1), \ldots, (\pi_n, g_n)) * ((\tau_1, h_1), \ldots, (\tau_n, h_n)) =
 ((\pi_{\tau_1}, g_{\tau_1} h_1), \ldots, (\pi_{\tau_n}, g_{\tau_n} h_n)).
\end{equation}
This description of $G\wr S_n$ is consistent with \cite{MR}. 

Given $u = ((\pi_1, g_1), \ldots, (\pi_n,g_n)) \in G \wr S_n$, let $\Co(u)$ denote the unique $G$-composition $((a_1, h_1), \ldots, (a_k,h_k))$ such that
\[ \begin{array}{ll}
\pi_1 < \pi_2 < \cdots < \pi_{a_1}, & g_1 = \cdots = g_{a_1} = h_1,\\
\pi_{a_1 + 1} < \cdots < \pi_{a_1 + a_2}, & g_{a_1+1} = \cdots = g_{a_1 + a_2} = h_2,\\
\vdots &  \\
\pi_{a_{1} + \cdots + a_{k-1}+1} < \cdots < \pi_n, & g_{a_1 + \cdots + a_{k-1}+1} = \cdots = g_n = h_k,
\end{array}
\]
and where $k$ is as small as possible. Thus, $\Co(u)$ keeps track of those values $i$ such that $\pi_i > \pi_{i+1}$ or $g_i \ne g_{i+1}$. 
For instance if $g, h$ are distinct elements in $G$, then
\[
\Co((3,g), (6, g), (4, g), (1, h), (2, h), (5, h), (8, g), (7, g))) = 
((2,g), (1, g), (3, h), (1, g), (1,g)).
\]
Let $\field[G\wr S_n]$ denote the group algebra of $G\wr S_n$.
For $\alpha\comp_G n$, define $Y_\alpha \in \field[G\wr S_n]$ by
\[ Y_\alpha = \sum_{u \in G\wr S_n : \Co(u) =\alpha} u. \]
Clearly the set $\{Y_{\alpha}~|~\alpha\comp_G n\}$ is linearly independent. Let
\[ \GDesAlg(G\wr S_n) = \text{$\field$-linear span of }  \{Y_\alpha~|~\alpha \comp_G n\}.\]
The following result is due to Mantaci and Reutenauer \cite[Theorem~6.9]{MR} (see also \cite{BH, BerH}):
\begin{theorem}\label{thm:MR}
If $G$ is abelian then $\GDesAlg(G \wr S_n)$ is a subalgebra of $\field[G\wr S_n]$.
\end{theorem}

A generalization of this theorem to arbitrary groups $G$ appears in \cite[Theorem~3.2]{NT}. These results can be deduced from our main theorem. First we will need to introduce another basis for $\GDesAlg(G \wr S_n)$. Consider the partial order on the set of $G$-compositions of $n$ generated by cover relations of the form
\[ ((a_1, g_1), \ldots, (a+b, g_i), \ldots, (a_k, g_k)) < ((a_1, g_1), \ldots, (a, g_i), (b, g_i), \ldots, (a_k, g_k)).\]
In other words $\alpha \le \beta$ if and only if $\beta$ is a color-preserving refinement of $\alpha$. For $\alpha \comp_G n$, let
\[ X_\alpha = \sum_{\beta \comp_G n:\beta \le \alpha} Y_\beta.\]
By M\"obius inversion,
\[ Y_\alpha = \sum_{\beta \comp_G n:\beta \le \alpha} (-1)^{\ell(\alpha) - \ell(\beta)} X_\beta.\]
Thus $(X_\alpha)_{\alpha \comp_G n}$ is a basis of $\GDesAlg(G\wr S_n)$.

\section{Proof of main result}

We restate and prove the main result announced in the Introduction. 

\bigskip\noindent
{\bf Theorem 1.} {\it The $\field$-module map $
f: (\field\Sigma_n^G)^{S_n} \to  \field[G\wr S_n]$ defined by
$ f(\sigma_\alpha) = X_\alpha$
is an injective anti-homomorphism of algebras.}

\bigskip
\begin{proof}
Let $\mathcal{C}$ be the set of ordered $G$-partitions of $[n]$ whose blocks are singletons. Note that $\mathcal{C}$ is a left ideal of $\Sigma_n^G$. To elaborate, given $P = ((B_1, h_1), (B_2, h_2), \ldots, (B_k, h_k)) \in \Sigma_n^G$ and $Q = ((\{\pi_1\}, g_1), (\{\pi_2\}, g_2), \ldots, (\{\pi_n\}, g_n))\in\mathcal{C}$, let $\tau\in S_n$ be the unique permutation such that 
\[
\begin{array}{ll}
B_1 = \{\pi_{\tau_1}, \pi_{\tau_2}, \ldots, \pi_{\tau_{a_1}}\} , &   \tau_1 < \cdots <\tau_{a_1},\\
B_2 = \{ \pi_{\tau_{a_1 + 1}}, \ldots, \pi_{\tau_{a_1 + a_2}}\},  &  \tau_{a_1 + 1} < \cdots < \tau_{a_1 + a_2}, \\
 \vdots & \\
B_k = \{\pi_{\tau_{a_1 + \cdots + a_{k-1} + 1}}, \ldots, \pi_{\tau_n}\}, & \tau_{a_1 + \cdots + a_{k-1} + 1} < \cdots < \tau_n.
\end{array}
\]
where $a_i = |B_i|$ for $i\in [k]$. Then it follows from the definition of multiplication in $\Sigma_n^G$ that
\begin{multline}\label{E:PQ}
PQ = ((\{ \pi_{\tau_1}\}, g_{\tau_1} h_1), (\{ \pi_{\tau_2}\}, g_{\tau_2} h_1), \ldots, (\{\pi_{\tau_{a_1}} \}, g_{\tau_{a_1}} h_1), \\
(\{\pi_{\tau_{a_1 + 1}} \}, g_{\tau_{a_1 + 1} }h_2), \ldots,  (\{\pi_{\tau_{a_1+a_2}}\}, g_{\tau_{a_1+a_2}}h_2), \\ \ldots, (\{ \pi_{\tau_{a_1 + \cdots + a_{k-1} + 1}}\}, g_{\tau_{a_1 + \cdots + a_{k-1} + 1}} h_k), \ldots, (\{ \pi_{\tau_n}\}, g_{\tau_n} h_k) ). 
\end{multline}

Consider the action of $(\field \Sigma_n^G)^{S_n}$ on the $\field$-module $\field \mathcal{C}$ by left multiplication. For any $\alpha = ((a_1, h_1), \ldots, (a_\ell, h_k)) \comp_G n$ and $((\{\pi_1\}, g_1), \ldots, (\{\pi_n\}, g_n)) \in \mathcal{C},$ by \eqref{E:PQ} we have
\begin{equation}\label{eq:LeftSemiMult}
 \sigma_\alpha ((\{\pi_1\},g_1), \ldots, (\{\pi_n\},g_n)) = 
\sum ((\{\pi_{\tau_1}\}, g_{\tau_1} i_1), \ldots, (\{\pi_{\tau_n}\}, g_{\tau_n} i_n))
\end{equation}
where the sum is over all $u = ((\tau_1, i_1), \ldots, (\tau_n, i_n)) \in G\wr S_n$ such that $\tau_1 < \cdots < \tau_{a_1},$ $\tau_{a_1+ 1} < \cdots < \tau_{a_1 + a_2},$ \ldots, $\tau_{a_1 + \cdots + a_{k-1} + 1} < \cdots < \tau_n$, and
$i_1 = \cdots = i_{a_1} = h_1,$ $i_{a_1 + 1} = \cdots = i_{a_1 + a_2} = h_2$,
\ldots, $i_{a_1 + \cdots + a_{k-1} + 1} = \cdots = i_n = h_k$. These conditions are equivalent to $\Co(u) \le \alpha$. 

Now identify $\mathcal{C}$ with the set $G\wr S_n$ so that if $v = ((\pi_1, g_1), \ldots, (\pi_n, g_n)) \in G\wr S_n$ then $v$ gets identified with $((\{\pi_1\}, g_1), \ldots, (\{\pi_n\}, g_n)) \in \mathcal{C}$. Let $I = ((1, e), (2, e), \ldots, (n, e))$, the identity element of $G\wr S_n$.
Comparing \eqref{eq:LeftSemiMult} with \eqref{E:desmult}, we get 
\[ \sigma_\alpha v = \sum_{u \in G\wr S_n :\Co(u) \le \alpha} v * u = v * (\sigma_\alpha \; I) \]
for any $v\in G\wr S_n$. In particular, $\sigma_\alpha I = X_\alpha$.

The map $f$ satisfies $f(\sigma_\alpha) = \sigma_\alpha I = X_\alpha$, and so $f(\sigma_\alpha \sigma_\beta) = \sigma_\alpha (\sigma_\beta I) = (\sigma_\beta I) * (\sigma_\alpha I) = X_\beta * X_\alpha$, completing the proof.
\end{proof}

Since the image of $f$ is $\GDesAlg(G\wr S_n)$, we have the following corollary:

\begin{corollary}
For any group $G$, $\GDesAlg(G\wr S_n)$ is a subalgebra of $\field[G\wr S_n]$ and is anti-isomorphic to $(\field\Sigma_n^G)^{S_n}$.
\end{corollary}

\section*{Acknowledgments} I thank Franco Saliola for invaluable comments on earlier versions of this manuscript, particularly for suggesting that I allow $G$ to be non-abelian. I also thank Ken Brown for helpful discussions, and Jean-Yves Thibon for bringing the paper \cite{NT} to my attention and providing instructive comments.

\end{document}